\documentclass[final]{ws-ccm}
\usepackage{amssymb,amsmath,mathrsfs}

\newcommand{\C}{{\mathbb C}}
\newcommand{\R}{{\mathbb R}}
\newcommand{\eps}{\varepsilon}

\begin{document}

\markboth{S.\ Secchi \& M.\ Squassina} {Schr\"odinger
Equations with Electromagnetic Field}

%
\catchline{}{}{}{}{}
%

\title{ON THE LOCATION OF SPIKES FOR THE SCHR\"ODINGER
EQUATION WITH ELECTROMAGNETIC FIELD}

\author{SIMONE SECCHI\footnote{The author was supported by the MIUR national research
project ``Variational Methods and Nonlinear Differential Equations''.}}

\address{Dipartimento di Matematica ``F.Enriques'',
Universit\`a di Milano \\
Via Saldini 50, I-20133 Milano, Italy \\
secchi@mat.unimi.it}

\author{MARCO SQUASSINA\footnote{The author was supported by the MIUR national research
project ``Variational and Topological Methods in the Study of Nonlinear Phenomena'' and
by the Istituto Nazionale di Alta Matematica ``F.Severi'' (INdAM)}}

\address{Dipartimento di Matematica ``F.Brioschi'',
Politecnico di Milano \\
Via Bonardi 9, I-20133 Milano, Italy \\
squassina@mate.polimi.it}

\maketitle

\begin{history}
\received{(Day Month Year)}
\revised{(Day Month Year)}
\end{history}

\begin{abstract}
We consider the standing wave solutions of the three dimensional semilinear
Schr\"odinger equation with competing potential functions $V$ and $K$ and
under the action of an external electromagnetic field $B$.
We establish some necessary conditions for a sequence
of such solutions to concentrate, in two different senses, around a
given point. In the particular but important
case of nonlinearities of power type, the spikes locate at
the critical points of a smooth ground energy map
independent of $B$.
\end{abstract}




\section{Introduction}
In this work we deal with the standing wave solutions
$$
\varphi(x,t)=e^{-\frac{iV_0}{\hbar}t}u(x), \qquad
x\in\R^3,\,\,t\in\R^+
$$
of the time-dependent Schr\"odinger equation with
electromagnetic field
$$
i\hbar\frac{\partial\varphi}{\partial
t}=\left(\textstyle\frac{\hbar}{i}\nabla
-A(x)\right)^2\varphi+W(x)\varphi-|\varphi|^{p-1}\varphi,
$$
where the Schr\"odinger operator is defined as
$$
\left(\textstyle\frac{\hbar}{i} \nabla -A\right)^2
:= -\hbar^2 \Delta-\textstyle\frac{2\hbar}{i}
\left\langle A \mid \nabla\right\rangle + |A|^2-
\textstyle\frac{\hbar}{i}\operatorname{div}A.
$$
Here $\hbar>0$ is the Planck constant, $p\in(1,5)$, and the
functions $W\colon\R^3\to\R$ and $A\colon\R^3\to\R^3$ are, respectively,
a scalar potential of the electric field $E=-\nabla W$ and
a vector potential for the external electromagnetic field $B={\rm curl}\,A$.
Now, the function $u\colon\R^3\to\C$ which appears
in $\varphi(x,t)$ satisfies, more generally, a time-independent
equation of the form
\begin{equation}
\label{schro}
\left(\textstyle\frac{\hbar}{i}\nabla-A(x)\right)^2u+V(x)u=K(x)f(|u|^2)u,
\end{equation}
where $V(x)=W(x)+V_0$, $K\colon\R^3\to\R$ is an additional potential function,
and $f\colon\R^+\to\R$ is a suitable nonlinearity. Quite recently, under
reasonable assumptions on $A$, $V$ and $K$, the study of the
existence of ground (bound) state solutions $u_{\hbar}$ to~\eqref{schro}
and the related investigation of the semi-classical limit
(the transition from Quantum to Classical Mechanics
as $\hbar\to 0$), has been tackled in various
contributions (see e.g.~\cite{az,chabro,cingo,cingosec,el,ku} for the case $A\neq 0$ and
\cite{blio,pinofelm,pinofelm3,pinofelm4,floerwein,oh1,rabinowitz} for the case $A=0$). More
precisely, it turns out that, if $z_0\in\R^3$ is a non-degenerate
critical point of the so called ground-energy
function $\Sigma_r\colon\R^3\to\R$ (see Definition~\ref{groundef}),
then for every $\hbar$ sufficiently
small~\eqref{schro} admits a least energy solution $u_{\hbar}$
concentrating near $z_0$. In the opposite direction, we are interested in
discussing some necessary conditions for the concentration of
a sequence of bound-state solutions to~\eqref{schro} in the neighborhood of a given point
$z_0$. In absence of the electromagnetic field, this problem has been
studied in various papers (see e.g.~\cite{abc,wang,wangzeng}),
mainly in the case where $f(u)$ is a power
of exponent $p$ (see also~\cite{g,pistoia}).
It turns out that, at least in this particular situation,
for the concentration to occur, $z_0$ has to be
a critical point for the $C^1$ ground-energy map (see \cite[Lemma 2.5]{wangzeng})
\begin{equation}
\label{explicit}
\Sigma_r(z)=\frac{V^{\frac{5-p}{2p-2}}(z)}{K^{\frac{2}{p-1}}(z)},
\qquad \text{for every $z\in\R^3$}.
\end{equation}
On the other hand, to our knowledge, for a more general nonlinearity
$f(u)$, the function $\Sigma_r(z)$ is locally Lipschitz continuous,
and its further smoothness properties seem to depend on the uniqueness
results for the limiting equation
\begin{equation}
\label{pb:frozen}
-\Delta u+V(z)u=K(z)f(|u|^2)u,
\end{equation}
where $z\in\R^3$ acts as a parameter. To overcome
this problem, recently, the authors have provided
in~\cite{secsqu} new necessary conditions involving generalized
derivatives of $\Sigma_r$ such as the Clarke subdifferential or even
weaker conditions, not requiring any regularity
of $\Sigma_r$ (see~Definition~\ref{gammapm}).
\par
Our purpose in this paper is to understand what happens
under the presence of an external electromagnetic vector potential
$A$, and to see whether $A$ may influence or not the
location of spikes for the solutions of~\eqref{schro}.
Actually, in general, this fact seems to depend on the notion of
concentration that one adopts. We consider at least two ways of
saying that a sequence $(u_\hbar)$ of bound-state solutions to~\eqref{schro}
is peaking around a given point
$z_0$.
The first one, the most intuitive, is a {\em pointwise concentration} and it is
precisely the one used in two papers by Wang and Zheng~\cite{wang,wangzeng}.
The second is a sort of {\em energetic concentration}
in terms of the functional associated with~\eqref{schro},
$$
J_{\hbar}(u)={\textstyle\frac{1}{2}}\!\int_{\R^3}|D^{\hbar}u|^2+V(x)|u|^2dx
-\int_{\R^3}K(x)F(|u|^2)dx,
$$
where $D^\hbar=\frac{\hbar}{i}\nabla-A(x)$. Precisely, we require that
$$
\lim_{\hbar\to 0}\hbar^{-3}J_{\hbar}(u_{\hbar})=\Sigma_r(z_0).
$$
As we prove in the main result, Theorem~\ref{main}, the vector potential $A$ might
affect the location of pointwise concentration points, whereas it does not influence
the energetic concentration points. In the particular but fairly significant case
where $f$ is a power nonlinearity, the above notions of concentration
coincide (see Proposition~\ref{comparison}), and it turns out that
the peaks locate at the classical critical points of the
smooth function~\eqref{explicit} independent of $A$, thus
rigorously confirming what conjectured in~\cite{cingosec}. In some sense,
from an heuristic point of view, $A$ tends
to lurk into the complex phase factor of the solutions.
We point out that, in the course of the proof
of Theorem~\ref{main}, we will derive an ad-hoc Pucci-Serrin
type identity for the complex-valued solutions to~\eqref{schro}
(cf.\ formula~\eqref{ps-prima}). Just for the sake of simplicity, we restrict the attention
to the physically relevant case of space-dimension $n=3$.
\bigskip


\begin{center}
\textbf{Notations}
\end{center}

\smallskip
\begin{enumerate}
\item $\Re w$ (resp.\ $\Im w$) stands for the real (resp.\ the
imaginary) part of $w\in\C$.
\item $i$ is the imaginary unit, namely $i^2 = -1$.
For $w\in\C$, we set $\bar w=\Re w - i \Im w$.
\item The gradient of a $C^1$ function $f:\R^3\to\R$ will be
denoted by $\nabla f$. The jacobian matrix of
a $C^1$ function $g:\R^3\to\R^3$ will be indicated by $g'$.
The directional derivatives of $f$ and $g$ along a vector $w$ will be indicated
by $\frac{\partial f}{\partial w}$ and $\frac{\partial g}{\partial w}$.
\item $\langle x\mid y\rangle$ denotes the standard scalar product in $\R^3$
of $x$ and $y$.
\end{enumerate}


\section{Problem Setting and Auxiliary Results}
\noindent
In this section, we collect a few preliminary definitions and results
that we need in order to state and prove our main achievement, Theorem~\ref{main}.
For the sake of simplicity, we rename the constant $\hbar$ into $\eps>0$.
We assume that the functions
$$
A\colon \R^3\to\R^3,\qquad
V\colon \R^3\to\R,\qquad
K \colon \R^3\to\R
$$
are all of class $C^1$, $K$ is positive and there exist $V_0>0$ and $K_0>0$ with
\begin{equation}
\label{VeK}
\inf_{x\in\R^3}V(x)=V_0
\qquad
\text{and}
\qquad
\sup_{x\in\R^3}K(x)=K_0.
\end{equation}
Moreover, the function $f:\R^+\to\R$ is of class $C^1$, increasing, $f(0)=0$ and
$$
\lim_{s\to\infty}\frac{f(s)}{s^{\frac{p-1}{2}}}=0
\quad\text{and}\quad
0<\vartheta F(s)\leq f(s)s
\quad\text{for some $p\in(1,5)$ and $\vartheta>2$},
$$
where $F(s)=\frac{1}{2}\int_0^s f(t)dt$ for $s\in\R^+$.
In order to formulate the problem in a suitable variational setting,
for every $\eps>0$, we introduce the (real) Hilbert space ${\mathcal H}^\eps_{A,V}$ defined
as the closure of $C^\infty_c(\R^3,\C)$ with respect to scalar product
$$
(u,v)_{{\mathcal H}^\eps_{A,V}}:=\Re\int_{\R^3}D^\eps u\overline{D^\eps v}+V(x)u\bar v\,dx,
\qquad\text{$D^\eps u=\textstyle\frac{\eps}{i}\nabla-A(x)$}.
$$
As remarked in~\cite{el}, ${\mathcal H}^\eps_{A,V}$ has in general no relationships with
$H^1(\R^3,\C)$. However, the following \textit{diamagnetic inequality}
is well known (see e.g.\ \cite{liblo})
\begin{equation}
\label{eq:dia}
\eps|\nabla |u|(x)| \leq |D^\eps u(x)|,
\qquad\text{for every $u\in {\mathcal H}^\eps_{A,V}$
and a.e.\ $x\in\R^3$},
\end{equation}
so that $|u|\in H^1(\R^3,\R)$ for any $u\in{\mathcal H}^\eps_{A,V}$.
Finally we recall that the Schr\"odinger operator is {\em gauge invariant}:
if we replace $A$ by $\tilde A=A+\nabla\chi$ for any
$\chi\in C^2(\R^3,\R)$, and we let $\tilde u=e^{\frac{i}{\eps}\chi}u$,
then ${\rm curl}\,\tilde A={\rm curl}\,A$ and
$$
\left(\textstyle\frac{\eps}{i} \nabla -\tilde A\right)\tilde u=
e^{\frac{i}{\eps}\chi}\left(\textstyle\frac{\eps}{i} \nabla -A\right)u,
$$
so that $\|\tilde u\|_{{\mathcal H}^\eps_{\tilde A,V}}\!\!=
\|u\|_{{\mathcal H}^\eps_{A,V}}.$
\vskip2pt
\noindent
Under the above assumptions, we give the following

\begin{definition}
We say that $(u_\eps)$ is a sequence of bound-state solutions to
\begin{equation}
\label{problem}
\tag{$S_\eps$}
\left(\textstyle\frac{\eps}{i}\nabla-A(x)\right)^2u+V(x)u=K(x)f(|u|^2)u
\end{equation}
if $u_\eps$ belongs to ${\mathcal H}^\eps_{A,V}$ for every $\eps>0$,
\begin{equation}
\label{boundL2}
\sup_{\eps>0}\eps^{-3}\|u_\eps\|^2_{{\mathcal H}^\eps_{A,V}}<\infty
\end{equation}
and $u_\eps$ satisfies~\eqref{problem} on $\R^3$ in weak sense.
\end{definition}

\subsection{The ground-energy functions}

Fixed $z\in\R^3$, we consider the functional
$$
I_{z}(u)={\textstyle\frac{1}{2}}\!\int_{\R^3}|\nabla u|^2+V(z)|u|^2dx
-\int_{\R^3}K(z)F(|u|^2)dx
$$
associated with the limiting equation~\eqref{pb:frozen}.
It is readily seen that $I_z$ is $C^1$
over both the spaces $H^1(\R^3,\R)$ and $H^1(\R^3,\C)$.

\begin{definition}
\label{groundef}
We define the real and the complex ground-state functions
$$
\Sigma_r \colon \R^3\to\R
\qquad
\text{and}
\qquad
\Sigma_c \colon \R^3\to\R
$$
by setting, for every $z\in\R^3$,
$$
\Sigma_r(z)=\min_{v\in\mathcal{N}_z}I_z(v)
\qquad
\text{and}
\qquad
\Sigma_c(z)=\min_{v\in\widetilde{\mathcal{N}}_z}I_z(v),
$$
where $\mathcal{N}_z$ (resp.\ $\widetilde{\mathcal{N}}_z$) are the real
(resp.\ the complex) Nehari manifolds,
\begin{align*}
\mathcal{N}_z &= \left\{ u\in H^1(\R^3,\R)\setminus \{0\} :\,\, I'_z (u)[u]=0 \right\}, \\
\widetilde{\mathcal{N}}_z &= \left\{ u\in H^1(\R^3,\C)
\setminus \{0\} :\,\, I'_z (u)[u]=0 \right\}.
\end{align*}
Here $I'_z(u)[v]$ stands for the directional derivative of $I_z$ at $u$ along $v$.
\end{definition}
\vskip2pt
\noindent

We denote by $S_r(z)$ the set of positive radial solutions up to translations
to~\eqref{pb:frozen} at the energy level $\Sigma_r(z)$.
As the next lemma claims, the map $\Sigma_r$ enjoys some
useful regularity properties (see~\cite{wangzeng}).

\begin{lemma}
\label{pro-sigma0}
The following facts hold:
\begin{description}
\item[(i)] $\Sigma_r$ is locally Lipschitz continuous;
\vskip3pt
\item[(ii)] the directional derivatives from the left and the
right of $\Sigma_r$ at every point $z\in\R^3$ along any $w\in\R^3$ exist and it holds
\begin{align*}
&\left(\frac{\partial\Sigma_r}{\partial w}\right)^{-}\!\!(z)=
\sup_{v\in S_r(z)}\left\langle\nabla_z I_z(v)\mid w\right\rangle, \\
&\left(\frac{\partial\Sigma_r}{\partial w}\right)^{+}\!\!(z)=
\inf_{v\in S_r(z)}\left\langle\nabla_z I_z(v)\mid w\right\rangle.
\end{align*}
Explicitly, we have
\begin{align*}
\left(\frac{\partial\Sigma_r}{\partial w}\right)^{-}\!\!(z)&=
\sup_{v\in S_r(z)}\Big[
\frac{\partial V}{\partial w}(z)
\int_{\R^3} \frac{|v|^{2}}{2}dx -
\frac{\partial K}{\partial w}(z)
\int_{\R^3} F(|v|^2)dx\Big], \\
\left(\frac{\partial\Sigma_r}{\partial w}\right)^{+}\!\!(z)&=
\inf_{v\in S_r(z)}\Big[
\frac{\partial V}{\partial w}(z)
\int_{\R^3} \frac{|v|^{2}}{2}dx-
\frac{\partial K}{\partial w}(z)
\int_{\R^3} F(|v|^2)dx\Big],
\end{align*}
for every $z,w\in\R^3$.
\end{description}
\end{lemma}

The next result will turn out to be pretty useful
along the proof of our main theorem. We stress that it
contains, as a particular case, Lemma 7 of \cite{ku}.
\begin{lemma}
\label{pro-sigma1}
The following facts hold:
\begin{description}
\item[(i)] $\Sigma_c(z)=\Sigma_r(z)$, for every $z\in\R^3$;
\item[(ii)] if $U_z:\R^3\to\C$ is a least energy solution
of problem~\eqref{pb:frozen}, then
$$
|\nabla |U_z|(x)|=|\nabla U_z(x)|
\qquad\text{and}\qquad
\Re\big(i\bar U_z(x)\nabla U_z(x)\big)=0,
$$
for a.e.\ $x\in\R^3$;
\item[(iii)] there exist $\omega\in\R$ and a real least
energy solution $u_z$ of problem~\eqref{pb:frozen} with
\begin{equation}
\label{c-r-rappr}
U_z(x)=e^{i\omega}u_z(x),\qquad
\text{for a.e.\ $x\in\R^3$}.
\end{equation}
\end{description}
\end{lemma}
\begin{proof}
Fix $z\in\R^3$. For the sake of convenience, we introduce the functionals
\begin{align*}
T(u) &= \int_{\R^3} |\nabla u|^2dx, \\
P_z(u) &= \int_{\R^3} \Big[ K(z)F(|u|^2) - {\textstyle\frac12} V(z) |u|^2 \Big]dx.
\end{align*}
Observe that $I_z (u) = \frac12 T(u) - P_z(u)$.
Consider the following minimization problems
\begin{align*}
\sigma_r(z)&=\min\big\{T(u): \text{$u\in H^1(\R^3,\R)$, $P_z(u)=1$}\big\}, \\
\sigma_c(z)&=\min\big\{T(u): \text{$u\in H^1(\R^3,\C)$, $P_z(u)=1$}\big\}.
\end{align*}
Note that, obviously, there holds $\sigma_c(z)\leq\sigma_r(z)$. If we denote
by $u_\star$ the Schwarz symmetric rearrangement (see e.g.~\cite{blio,liblo})
of the positive real valued function $|u|\in H^1(\R^3,\R)$, then,
Cavalieri's principle yields
$$
\int_{\R^3} F(|u_\star|^2)dx =\int_{\R^3} F(|u|^2)dx
\qquad
\text{and}
\qquad
\int_{\R^3} |u_\star|^2dx =\int_{\R^3} |u|^2dx,
$$
which entails $P_z(u_\star)=P_z(u)$. Moreover, by the Polya-Szeg\"o inequality, we have
$$
T(u_\star)=\int_{\R^3} |\nabla u_\star|^2dx \leq \int_{\R^3} |\nabla |u||^2dx
\leq \int_{\R^3} |\nabla u|^2dx=T(u),
$$
where the second inequality follows by~\eqref{eq:dia} with $A=0$ and $\eps=1$.
Therefore, one can compute $\sigma_c(z)$ by minimizing over
the subclass of positive, radially symmetric and radially
decreasing functions $u\in H^1(\R^3,\R)$. As a consequence,
$\sigma_r(z)\leq\sigma_c(z)$. In conclusion, $\sigma_r(z)=\sigma_c(z)$.
Observe now that
\begin{align*}
\Sigma_r(z)&=\min\big\{I_z (u): \text{$u\in H^1(\R^3,\R)\setminus\{0\}$
is a solution to~\eqref{pb:frozen}}\big\}, \\
\Sigma_c(z)&=\min\big\{I_z (u): \text{$u\in H^1(\R^3,\C)\setminus\{0\}$
is a solution to~\eqref{pb:frozen}}\big\}.
\end{align*}
The above equations hold since any nontrivial real (resp.\ complex) solution of
\eqref{pb:frozen} belongs to $\mathcal{N}_z$ (resp.\ $\widetilde{\mathcal{N}_z}$) and,
conversely, any solution of $\Sigma_r(z)$ (resp.\ $\Sigma_c(z)$) produces a nontrivial
solution of \eqref{pb:frozen}. Moreover, it follows from an
easy adaptation of~\cite[Theorem 3, p.331]{blio} that
$\Sigma_r(z)=\sigma_r(z)$ as well as $\Sigma_c(z)=\sigma_c(z)$. In conclusion,
$$
\Sigma_r(z)=\sigma_r(z)=\sigma_c(z)=\Sigma_c(z),
$$
which proves (i). To prove (ii), let $U_z:\R^3\to\C$ be a least energy solution
to problem~\eqref{pb:frozen}. There holds
$|\nabla |U_z||\leq |\nabla U_z|$. Assume by contradiction that
$$
{\mathcal L}^3\big(\left\{x\in\R^3:
|\nabla |U_z|(x)|<|\nabla U_z(x)|\right\}\big)>0,
$$
where ${\mathcal L}^3$ is the Lebesgue measure in $\R^3$.
Then we get $P_z(|U_z|)=P_z(U_z)$ and
$$
\sigma_r(z)\leq\int_{\R^3} |\nabla |U_z||^2dx
<\int_{\R^3} |\nabla U_z|^2dx=\sigma_c(z),
$$
which is a contradiction. The second assertion in (ii) follows
by a direct computation. Indeed, a.e.\ in $\R^3$, we have
$$
|\nabla |U_z||=|\nabla U_z|
\qquad
\text{if and only if}
\qquad
\Re U_z\nabla (\Im U_z)=\Im U_z\nabla(\Re U_z).
$$
If this last condition holds, in turn, a.e.\ in $\R^3$ we have
$$
{\bar U}_z\nabla U_z=\Re U_z\nabla (\Re U_z)+
\Im U_z\nabla (\Im U_z),
$$
which implies the desired assertion. Finally, the
representation formula of (iii) is an immediate consequence
of (ii), since one obtains $U_z=e^{i\omega}|U_z|$ for some $\omega\in\R$.
\end{proof}

\subsection{Generalized gradients}
\label{gengrad}
Assume that $f:\R^3\to\R$ is a locally Lipschitz continuous function.
For the reader convenience, we recall that the Clarke
subdifferential (or generalized gradient) of $f$
at a point $z$ (cf.\ \cite{clarke}) is defined as
$$
\partial_C f(z)=\Big\{\eta\in\R^3:\,\,
f^0(z,w)\geq \left\langle\eta\mid w\right\rangle,\,\,\,\text{for every $w\in\R^3$}\Big\},
$$
where $f^0(z,w)$ is the Clarke derivative of
$f$ at $z$ along the direction $w$, defined as
$$
f^0(z;w)=\limsup_{\substack{\xi \to z
\\ \lambda \to 0^+}}\frac{f(\xi+\lambda w)-f(\xi)}{\lambda}.
$$
From~\cite[Proposition 2.3.1]{clarke} we learn that
$\partial_C f(z)$ is nonempty, convex and
\begin{equation}
\label{oppositesub}
\partial_C (-f)(z)=-\partial_C f(z),\qquad\text{for every $z\in\R^3$}.
\end{equation}
\noindent
In light of (i) in Lemma~\ref{pro-sigma0}, we are allowed to give the following
\begin{definition}
We denote by $\mathfrak{S}\subset\R^3$ the set of critical points of the function $\Sigma_r$
in the sense of the Clarke subdifferential, namely
$$
\mathfrak{S}:=\big\{z\in\R^3:\,\, 0\in\partial_C\Sigma_r(z)\big\}.
$$
\end{definition}

\noindent
Now, for $z\in\R^3$, we consider the gauge invariant functional
$J_z\colon H^1(\R^3,\C)\to\R$
$$
J_{z}(u) = {\textstyle\frac{1}{2}}\int_{\R^3} \left| \left( {\textstyle\frac{1}{i}}\nabla
-A(z) \right) u \right|^2 + V(z)|u|^2 dx-\int_{\R^3}K(z)F(|u|^2)dx,
$$
associated with the limiting equation
$$
\left({\textstyle\frac{1}{i}}\nabla-A(z)\right)^2u+V(z)u=K(z)f(|u|^2)u.
$$
We denote by $G_c(z)$ the set of the nontrivial
solutions $v:\R^3\to\C$, up to translations, of the
above limiting problem with bounded, but not necessarily least, energy.
Moreover, we introduce the linear map $\Upsilon_z:\R^3\to\R$, defined as
$$
\Upsilon_z(x):=\sum_{j=1}^3A_j(z)x_j,
\qquad\text{for every $x\in\R^3$}.
$$
Apparently, for every $z\in\R^3$, there holds $\nabla\Upsilon_z(x)=A(z)$.
It is readily seen that for every $v\in G_c(z)$ we can write $v=e^{i\Upsilon_z}U_z$, where
$U_z$ is a (possibly complex) solution to problem~\eqref{pb:frozen}.

\begin{definition}
\label{gammapm}
Let $z\in\R^3$. For every $w\in\R^3$ we define $\Gamma^-_z(w)$
and $\Gamma^+_z(w)$ by
$$
\Gamma^-_z(w):=\sup_{v\in G_c(z)}\left\langle\nabla_z J_z(v)\mid w\right\rangle
\quad\text{and}\quad
\Gamma^+_z(w):=\inf_{v\in G_c(z)}\left\langle\nabla_z J_z(v)\mid w\right\rangle,
$$
where $\nabla_z$ is the gradient with respect to $z$.
Explicitly, for every $w\in\R^3$,
\begin{align*}
\Gamma_z^-(w)&=
\sup_{\substack{v=e^{i\Upsilon_z}U_z \\ v\in G_c(z)}}\Big[
\Big\langle\frac{\partial A}{\partial w}(z)\mid\int_{\R^3} \Re(i\bar U_z
\nabla U_z)dx\Big\rangle  \\
&\qquad +\frac{\partial V}{\partial w}(z)
\int_{\R^3} \frac{|U_z|^{2}}{2}dx-
\frac{\partial K}{\partial w}(z)
\int_{\R^3} F(|U_z|^2)dx\Big], \\
\Gamma_z^+(w)&=
-\inf_{\substack{v=e^{i\Upsilon_z}U_z \\ v\in G_c(z)}}\Big[
\Big\langle\frac{\partial A}{\partial w}(z)\mid\int_{\R^3} \Re(i\bar U_z
\nabla U_z)dx\Big\rangle  \\
&\qquad +\frac{\partial V}{\partial w}(z)
\int_{\R^3} \frac{|U_z|^{2}}{2}dx-
\frac{\partial K}{\partial w}(z)
\int_{\R^3} F(|U_z|^2)dx\Big].
\end{align*}
Notice that
$$
\partial\Gamma_z^\pm(0)=\Big\{\eta\in\R^3:\,\,
\text{$\Gamma_z^\pm(w)\geq \left\langle\eta\mid w\right\rangle$,\,\,
for every $w\in\R^3$}\Big\},
$$
where $\partial\Gamma_z^\pm(0)$ is the subdifferential of the convex
function $\Gamma_z^\pm$ at zero. We set
$$
\mathfrak{S}^*:=\big\{z\in\R^3:\,\,0\in\partial\Gamma_z^-(0)\cap\partial\Gamma_z^+(0)\big\}
$$
and we say that $\mathfrak{S}^*$ is the set of
weak-concentration points for problem~\eqref{problem}.
\end{definition}

\subsection{Concentration of bound-state solutions}
We now introduce two (gauge invariant) notions of concentration
for a sequence of bound-states solutions of~\eqref{problem}
around a given point.

\begin{definition}
\label{conc-defs}
Let $z_0\in\R^3$ and assume that $(u_{\eps_h})\subset {\mathcal H}^{\eps_h}_{A,V}$
is a sequence of bound-state solutions to problem~\eqref{problem}. We say that
\begin{description}
\item[(i)] $z_0$ is a {\em concentration point}
for $(u_{\eps_h})$ if $|u_{\eps_h}(z_0)|\geq\varrho>0$ and
for every $\eta>0$ there exist $\rho>0$ and $h_0\geq 1$ such that
$$
|u_{\eps_h}(x)|\leq \eta,\qquad
\text{for every $h\geq h_0$ and $|x-z_0|\geq\eps_h\rho$}.
$$
The set of such points will be denoted by $\mathscr{C}\subset\R^3$;
\vskip2pt
\item[(ii)] $z_0$ is an {\em energy-concentration point} if
$$
\lim_{h\to \infty} \eps_h^{-3} J_{\eps_h}(u_{\eps_h})=\Sigma_r(z_0).
$$
The set of such points will be denoted by $\mathscr{E}\subset\R^3$.
\end{description}
\end{definition}

\noindent
For instance, if $K\equiv 1$, $f$ is a power, $z_0$ is a minimum point of $V$
and $(u_{\eps_h})$ is a sequence of least-energy solutions
to~\eqref{problem}, then $z_0\in\mathscr{E}\neq\emptyset$ (cf.\ ~\cite[Lemma 3]{ku}).

\vskip2pt
\noindent
Next we see that in the case of power nonlinearities
\begin{equation}
\label{fpowernon}
f(u)=\lambda u^{\frac{p-1}{2}}\qquad\text{for some $p\in(1,5)$ and $\lambda>0$},
\end{equation}
the above notions (i) and (ii) coincide.

\begin{proposition}
\label{comparison}
Let $f$ be as in~\eqref{fpowernon}. Then $\mathscr{E}=\mathscr{C}$.
\end{proposition}
\begin{proof}
Let $z_0 \in\mathscr{E}$ and consider $v_h(x)=u_{\eps_h}(z_0+\eps_hx)$.
Then $(|v_h|)$ converges to some $\tilde{v}\geq 0$ weakly in $H^1(\R^3,\R)$
and strongly in $L_{\rm loc}^q(\R^3,\R)$ for $2\leq q<6$
(see Step I in the proof of Theorem~\ref{main}).
By Kato's inequality~\cite[Theorem X.33]{rs}, we get
$$
\int_{\R^3} K(z_0+\eps_h x) |v_h|^p\tilde v \,dx
\geq \int_{\R^3} \nabla |v_h|\nabla\tilde v+V(z_0+\eps_h x)|v_h|\tilde v \,dx
$$
which, as $h\to\infty$, yields,
$$
\int_{\R^3} K(z_0) |\tilde v|^{p+1}dx
\geq \int_{\R^3} |\nabla\tilde v|^2+V(z_0)|\tilde v|^2dx.
$$
Therefore, there exists $\vartheta\in(0,1]$ such that
$\vartheta \tilde{v}\in \mathcal{N}_{z_0}$. As a consequence,
\begin{align*}
\Sigma_r(z_0) &\leq \vartheta^2 \left(\textstyle{ \frac12 - \frac1{p+1}} \right) \int_{\R^3}
|\nabla \tilde{v}|^2 + V(z_0)|\tilde{v}|^2 dx \\
&\leq \left(\textstyle{ \frac12 - \frac1{p+1}} \right) \liminf_{h\to\infty}
\int_{\R^3} | \nabla |v_h||^2 + V(z_0+\eps_h x) |v_h|^2 dx \\
&\leq \left(\textstyle{ \frac12 - \frac1{p+1}} \right) \liminf_{h\to\infty}
\int_{\R^3} \left| \left( \textstyle{\frac{1}{i}}\nabla - A(z_0+\eps_h x) \right) v_h
\right|^2 + V(z_0+\eps_h x)|v_h|^2 dx \\
\noalign{\vskip2pt}
&\leq \liminf_{h\to\infty} \eps_h^{-3} J_{\eps_h}(u_{\eps_h})=\Sigma_r(z_0),
\end{align*}
where we have used the diamagnetic inequality~\eqref{eq:dia} with $\eps=1$.
Hence we get $\vartheta=1$, which gives at once $\tilde{v}\in \mathcal{N}_{z_0}$. Then,
\begin{align*}
\int_{\R^3} |\nabla \tilde{v}|^2 + V(z_0) |\tilde{v}|^2dx &\leq
\left(\textstyle{ \frac12 - \frac1{p+1}} \right)^{-1}\liminf_{h\to\infty} \eps_h^{-3}
J_{\eps_h}(u_{\eps_h}) \\
&= \left(\textstyle{ \frac12 - \frac1{p+1}} \right)^{-1} \Sigma_r (z_0)
\leq\int_{\R^3} |\nabla \tilde{v}|^2 + V(z_0) |\tilde{v}|^2dx.
\end{align*}
This implies that $|v_h| \to \tilde{v}$
strongly in $H^1(\R^3,\R)$. Repeating the arguments
in the proof of \cite[Lemma 5]{ku} we conclude that $z_0 \in \mathscr{C}$
(the concentration occurs exponentially fast, see Step II of
the proof of Theorem~\ref{main}). This proves that $\mathscr{E}\subset\mathscr{C}$.
The converse inclusion follows by the uniqueness
of solutions (up to translations) to problem~\eqref{pb:frozen}. Indeed,
if $z_0 \in \mathscr{C}$, the sequence $\eps_h^{-3} J_{\eps_h}(u_{\eps_h})$
converges to $J_{z_0}(v_0)$ being $v_0$ an element of the family
$$
\{e^{i\Upsilon_{z_0}(x)+i\omega}\phi_0(x)\}_{\omega\in\R},
$$
where $\phi_0$ is the unique solution to~\eqref{pb:frozen} up
to translations (cf.\ \cite[Lemma 7]{ku}). In particular,
there holds $J_{z_0}(v_0)=I_{z_0}(\phi_0)=\Sigma_r(z_0)$, that is $z_0\in\mathscr{E}$,
concluding the proof. For similar considerations in the case $A=0$, see
e.g.~Lemma 4.2 in~\cite{gm}.
\end{proof}

\noindent
We are naturally lead to consider the following question
(see also Remark~\ref{rem-gauge}).
\begin{question}
\label{op-quest}
When $f(u)$ does {\em not} satisfy~\eqref{fpowernon}, is it
still true that $\mathscr{E}=\mathscr{C}$?
\end{question}

\section{The Main Result}

\noindent
For every $p\in(1,5)$, let us set
$$
\mathfrak{S}_p:=\big\{z\in\R^3:(5-p)K(z)\nabla V(z)=4V(z)\nabla K(z)\big\}.
$$
We now come to the main result of the paper.

\begin{theorem}
\label{main}
Assume that there exist $C\geq 0$ and $\gamma>0$ such that, for $|x|$ large,
\begin{equation}
\label{growthexpo}
\left|A'(x)\right|\leq Ce^{\gamma|x|},
\quad
\left|\nabla V(x)\right|\leq Ce^{\gamma|x|},
\quad
\left|\nabla K(x)\right|\leq Ce^{\gamma|x|}.
\end{equation}
Let $(u_{\eps_h})\subset {\mathcal H}^{\eps_h}_{A,V}$ be a sequence of
bound-state solutions to~\eqref{problem}. Then,
$$
\mathscr{C}\subset\mathfrak{S}^*
\qquad
\text{and}
\qquad
\mathscr{E}\subset\mathfrak{S}.
$$
If in addition $f$ satisfies~\eqref{fpowernon}, then we have
$$
\mathscr{C}=\mathscr{E}\subset \mathfrak{S}=\mathfrak{S}_p.
$$
\end{theorem}

\begin{proof}
Let $z_0\in\mathscr{C}$ and set $v_h(x)=u_{\eps_h}(z_0+\eps_hx)$
for every $h\geq 1$ and $x\in\R^3$. Then,
the sequence $(v_h)$ satisfies the rescaled equation
\begin{align}
\label{eq-h}
-\Delta v_h &- {\textstyle\frac{2}{i}}\left\langle A(z_0+\eps_hx)
\mid\nabla v_h\right\rangle -{\textstyle\frac{\eps_h}{i}}{\rm div}A(z_0+\eps_hx)v_h  \notag \\
&+|A(z_0+\eps_hx)|^2v_h+V(z_0+\eps_hx)v_h=K(z_0+\eps_hx)f(|v_h|^2)v_h.
\end{align}
We shall divide the proof into five steps.
\vskip2pt
\noindent
\textbf{Step I.} Up to a subsequence, $(v_h)$ converges in some
H\"{o}lder space $C_{\rm loc}^{2,\alpha}(\R^3)$ to the function
$v_0(x)=e^{i\Upsilon_{z_0} (x)} U_{z_0}(x)$, where $U_{z_0}:\R^3\to\C$ is a
solution to the equation
\begin{equation}
\label{limU}
-\Delta U_{z_0} + V(z_0) U_{z_0} = K(z_0) f(|U_{z_0}|^2) U_{z_0}.
\end{equation}
By the assumption on $(u_{\eps_h})$, the sequence
$(v_h)$ is bounded in ${\mathcal H}^{1}_{A,V}$,
and the diamagnetic inequality \eqref{eq:dia} immediately implies
that $(|v_h|)$ is bounded in $H^1(\R^3,\R)$.
Therefore, up to a subsequence, it converges weakly in $H^1(\R^3,\R)$ and
locally strongly in any $L^q(\R^3,\R)$ with $q<6$ towards
a positive function $v_*$. Moreover, for each
compact subset $\Lambda\subset\R^3$, by the continuity of $A$,
$(v_h)$ is also bounded in $H^1(\Lambda,\C)$.
We may now use the subsolution estimate (see e.g.\ \cite[Theorem 8.17]{gt}) to get
that $(v_h)$ is also bounded in $L^\infty_{{\rm loc}}(\R^3)$
and hence in $C^{2,\alpha}_{{\rm loc}}(\R^3)$, via Schauder' estimates. By combining this
fact with the results of \cite{lu}, up to a
subsequence, $v_h$ converges to $v_0$ in $C^{2,\alpha}_{{\rm loc}}(\R^3)$ and
furthermore $v_0\not\equiv 0$, since $|v_h(0)|=|u_{\eps_h}(z_0)|\geq \varrho>0$.
By continuity, the limit $v_0$ satisfies the limiting equation
\begin{equation}
\label{3}
-\Delta v_0 - {\textstyle\frac{2}{i}} \left\langle A(z_0)
\mid \nabla v_0\right\rangle + |A(z_0)|^2 v_0 + V(z_0) v_0
= K(z_0) f(|v_0|^2)v_0.
\end{equation}
If we define $U_{z_0}: x\in \R^3 \mapsto e^{-i \Upsilon_{z_0} (x)} v_0(x)$, then
$U_{z_0}$ satisfies~\eqref{limU}.
\vskip2pt
\noindent
\textbf{Step II.} There exist two positive constants
$R_*$ and $C_*$ such that
\begin{equation}
\label{expdecay}
|v_h(x)| \leq C_* e^{-\sqrt{\frac{V_0}{2}} |x|},
\qquad\text{for every $|x|\geq R_*$ and $h\geq 1$},
\end{equation}
where $V_0$ is defined in \eqref{VeK}. Since $z_0\in\mathscr{C}$, we have $v_h(x)\to 0$
as $|x| \to \infty$, uniformly with respect to $h\geq 1$.
Hence, for any $\eta>0$, we can find a radius $R_\eta>0$ such
that $|v_h(x)|<\eta$ whenever $|x|>R_\eta$ and $h\geq 1$. Therefore,
exploiting again Kato's inequality
$$
\Delta |v_h|\geq
\Re(\bar v_h|v_h|^{-1}(\nabla-iA)^2v_h)\qquad\text{(in distributional sense)},
$$
and taking into account that $f$ is increasing, there holds
$$
\Delta |v_h| \geq V(z_0 +\eps_h x) |v_h| - K(z_0+\eps_h x) f(|v_h|^2) |v_h|
\geq [V_0- K_0 f(\eta^2)] |v_h|
$$
in the sense of distributions on $\{|x|>R_\eta\}$, where $K_0>0$
is as in \eqref{VeK}. Let $\Gamma_0$ be a
fundamental solution for $-\Delta + c_\eta$, where $c_\eta=V_0-K_0f(\eta^2)$.
We can choose $\Gamma_0$ so that $|v_h (x)|
\leq [V_0-K_0f(\eta^2)] \Gamma_0 (x)$ holds for $|x|=R_\eta$.
Then, if $w=|v_h|-[V_0-K_0f(\eta^2)]\Gamma_0$, there holds
\begin{align*}
\Delta w &= \Delta |v_h| - [V_0-K_0f(\eta^2)] \Delta \Gamma_0 \\
&\geq [V_0-K_0 f(\eta^2)]|v_h| - [V_0-K_0f(\eta^2)]^2 \Gamma_0 \\
\noalign{\vskip1pt}
&= [V_0-K_0f(\eta^2)] w
\end{align*}
in distributional sense over $\{|x|>R_\eta\}$. Then, by the maximum principle,
$w(x) \leq 0$ for every $|x| \geq R_\eta$. Since,
as known, $\Gamma_0$ decays exponentially at the rate $\sqrt{c_\eta}$,
fixing $\eta=\eta_*$ so small that $f(\eta_*^2)\leq V_0/2K_0$,
we can find constants $R_*>0$ and $c>0$ such that $\Gamma_0(x)
\leq c\exp\{-\sqrt{V_0/2}|x|\}$ for $|x| \geq R_*$,
which yields the desired conclusion.
\vskip2pt
\noindent
\textbf{Step III.} For every $h\geq 1$, the
following identity holds
\begin{align}
\label{ps-prima}
\int_{\R^3}\bigg[
& \Big\langle\frac{\partial A}{\partial x_k}(z_0+\eps_h x)\mid
A(z_0+\eps_h x)\Big\rangle|v_h|^2-\Re\Big\langle{\textstyle\frac{1}{i}}\nabla v_h \mid
\frac{\partial A}{\partial x_k}(z_0+\eps_h x)\bar v_h\Big\rangle \notag \\
&+ \frac{\partial V}{\partial x_k}
(z_0+\eps_hx)\frac{|v_h|^2}{2}-\frac{\partial K}{\partial x_k}
(z_0+\eps_hx) F(|v_h|^2)\bigg]dx =0.
\end{align}
Rigorously, we cannot directly apply
the Pucci-Serrin variational identity~\cite{ps}, since the
solutions to equation~\eqref{eq-h} are complex-valued. For we are not aware of any explicit
reference to cite for the identity we need, we will derive~\eqref{ps-prima} directly (see
also~\cite{colin}).
Throughout the rest of this step only, we use
the less cumbersome notation $x\cdot y$ in place of $\langle x\mid y\rangle$ to indicate the
standard scalar product in $\R^3$.

\noindent
First of all, let us observe that, for every $h\geq 1$,
$$
|\nabla v_h| \leq |D^1 v_h| + |A(z_0+\eps_h x)| |v_h|.
$$
Hence, taking into account Step II and the
bounds~\eqref{growthexpo} and~\eqref{boundL2}, we get
\begin{align}
\label{boungrad}
& \|\nabla v_h\|_{L^2(\R^3)}\leq\| D^1 v_h \|_{L^2(\R^3)}
+\| A(z_0+\eps_h x) v_h \|_{L^2(\R^3)} \notag \\
&\leq \| D^1 v_h \|_{L^2(\R^3)}
+|A(z_0)|\|v_h\|_{L^2(\R^3)}+c\|e^{\gamma\eps_h |x|}|x|v_h \|_{L^2(\R^3)} \leq c,
\end{align}
for all $h\geq 1$ and some $c>0$.
Let $\delta>0$ and consider
the cut-off function $\psi_\delta=\psi(\delta x)$, where
$\psi\in C_c^1(\R^3)$ is such that $\psi(x)=1$ for $|x|\leq 1$ and $\psi(x)=0$
for $|x|\geq 2$. If $\boldsymbol{e}_k$
denotes the $k$-th vector of the canonical base in $\R^3$, we test equation~\eqref{eq-h}
with the function $\psi_\delta\boldsymbol{e}_k\cdot\overline{\nabla v_h}$ and we
take the real part. Firstly, we have
\begin{align*}
\Re\int_{\R^3}\nabla v_h\cdot\nabla[\psi_\delta
\boldsymbol{e}_k\cdot\overline{\nabla v_h}]dx=
\Re\int_{\R^3}\nabla v_h\cdot\nabla\psi_\delta
\boldsymbol{e}_k\cdot\overline{\nabla v_h}dx
-\int_{\R^3}\nabla\psi_\delta\cdot\boldsymbol{e}_k
\frac{|\nabla v_h|^2}{2}dx.
\end{align*}
As a consequence, by virtue of~\eqref{boungrad},
the Dominated Convergence Theorem yields
$$
\lim_{\delta\to 0}\Re\int_{\R^3}\nabla v_h\cdot\nabla[\psi_\delta
\boldsymbol{e}_k\cdot\overline{\nabla v_h}]dx=0.
$$
Now, we have
\begin{align*}
&\Re\int_{\R^3}K(z_0+\eps_h x)f(|v_h|^2)v_h\psi_\delta
\boldsymbol{e}_k\cdot\overline{\nabla v_h}dx=\Re\int_{\R^3}K(z_0+\eps_h x)\psi_\delta
\boldsymbol{e}_k\cdot \nabla F(|v_h|^2)dx \\
&=-\eps_h\int_{\R^3}\frac{\partial K}{\partial x_k}(z_0+\eps_h x)\psi_\delta
F(|v_h|^2)dx-\int_{\R^3}K(z_0+\eps_h x)\frac{\partial\psi_\delta}{\partial x_k}F(|v_h|^2)dx.
\end{align*}
Hence, in light of~\eqref{growthexpo},~\eqref{expdecay} and~\eqref{boungrad},
by the Dominated Convergence Theorem we have
$$
\lim_{\delta\to 0}\Re\int_{\R^3}K(z_0+\eps_h x)f(|v_h|^2)v_h\psi_\delta
\boldsymbol{e}_k\cdot\overline{\nabla v_h}dx=-\eps_h\int_{\R^3}
\frac{\partial K}{\partial x_k}(z_0+\eps_h x)F(|v_h|^2)dx.
$$
In a similar fashion, there hold
$$
\lim_{\delta\to 0}\Re\int_{\R^3}V(z_0+\eps_h x)v_h\psi_\delta
\boldsymbol{e}_k\cdot\overline{\nabla v_h}dx=-\eps_h\int_{\R^3}
\frac{\partial V}{\partial x_k}(z_0+\eps_h x)\frac{|v_h|^2}{2}dx,
$$
$$
\lim_{\delta\to 0}\Re\int_{\R^3}|A(z_0+\eps_h x)|^2v_h\psi_\delta
\boldsymbol{e}_k\cdot\overline{\nabla v_h}dx=-\eps_h\int_{\R^3}
A(z_0+\eps_h x)\cdot\frac{\partial A}{\partial x_k}(z_0+\eps_h x)|v_h|^2dx.
$$
Finally, we have
$$
J(\delta)=\Re\int_{\R^3}{\textstyle\frac{2}{i}}A(z_0+\eps_h x)\cdot\nabla v_h\psi_\delta
\boldsymbol{e}_k\cdot\overline{\nabla v_h}dx=J_1(\delta)+J_2(\delta)+J_3(\delta),
$$
where we have set
\begin{align*}
J_1(\delta)&=-\eps_h\Re\sum_{m=1}^3\int_{\R^3}{\textstyle\frac{2}{i}}
\frac{\partial A_m}{\partial x_k}(z_0+\eps_h x)
\psi_\delta\frac{\partial v_h}{\partial x_m}\bar v_hdx, \\
J_2(\delta)&=-\Re\sum_{m=1}^3\int_{\R^3}{\textstyle\frac{2}{i}}A_m(z_0+\eps_h x)
\frac{\partial\psi_\delta}{\partial x_k}\frac{\partial v_h}{\partial x_m}\bar v_hdx, \\
J_3(\delta)&=-\Re\sum_{m=1}^3\int_{\R^3}{\textstyle\frac{2}{i}}A_m(z_0+\eps_h x)
\psi_\delta\frac{\partial^2 v_h}{\partial x_k\partial x_m}
\bar v_hdx.
\end{align*}
After a few computations, one shows
that $J_2(\delta)\to 0$ as $\delta\to 0$ and
$$
J_3(\delta)=-\Re\int_{\R^3}
{\textstyle\frac{2\eps_h}{i}}{\rm div}A(z_0+\eps_h x)v_h\psi_\delta
\boldsymbol{e}_k\cdot\overline{\nabla v_h}dx-J(\delta)+\Theta(\delta),
$$
with $\Theta(\delta)\to 0$ as $\delta\to 0$. Furthermore,
again by~\eqref{growthexpo},~\eqref{expdecay}
and~\eqref{boungrad}
$$
\lim_{\delta\to 0}J_1(\delta)=-\eps_h\Re\int_{\R^3}{\textstyle\frac{2}{i}}\nabla v_h \cdot
\frac{\partial A}{\partial x_k}(z_0+\eps_h x)\bar v_h dx.
$$
Therefore, we obtain
$$
\lim_{\delta\to 0}J(\delta)=-
\Re\int_{\R^3}
{\textstyle\frac{\eps_h}{i}}{\rm div}A(z_0+\eps_h x)v_h
\boldsymbol{e}_k\cdot\overline{\nabla v_h}dx
-\eps_h\Re\int_{\R^3}{\textstyle\frac{1}{i}}\nabla v_h \cdot
\frac{\partial A}{\partial x_k}(z_0+\eps_h x)\bar v_h dx.
$$
Adding the above identities immediately yields~\eqref{ps-prima}.
\vskip2pt
\noindent
\textbf{Step IV.} We apply the Dominated Convergence Theorem to take the limit
as $h\to\infty$ into identity~\eqref{ps-prima}. The only troublesome term is
$$
\Re\Big\langle{\textstyle\frac{1}{i}} \nabla v_h \mid
\frac{\partial A}{\partial x_k}(z_0+\eps_h x)\bar v_h\Big\rangle,
$$
since we apparently have no control on the decay of $\nabla v_h$.
Taking into account~\eqref{boungrad} and recalling that
$\nabla v_h(x)\to\nabla v_0(x)$ for all $x\in\R^3$,
up to a subsequence, we have
\begin{equation}
\label{debole}
\nabla v_h\rightharpoonup \nabla v_0,
\qquad\text{weakly in $L^2(\R^3)$}.
\end{equation}
On the other hand, by virtue of Step II, there exist
$R_*>0$ and $c>0$ such that
$$
\left|\frac{\partial A}{\partial x_k}(z_0+\eps_h x)\bar v_h\right|
\leq ce^{-\big(\sqrt{\frac{V_0}{2}}-\gamma\eps_h\big) |x|},
\qquad\text{for every $|x|\geq R_*$}.
$$
Consequently, since $\bar v_h(x)\to \bar v_0(x)$ for all
$x\in\R^3$ and $A\in C^1(\R^3)$, there holds
\begin{equation}
\label{forte}
\frac{\partial A}{\partial x_k}(z_0+\eps_h x)\bar v_h\to
\frac{\partial A}{\partial x_k}(z_0)\bar v_0,
\qquad\text{strongly in $L^2(\R^3)$}.
\end{equation}
Thus, by combining~\eqref{debole} and~\eqref{forte}, for each $k$, we immediately get
$$
\lim_{h\to\infty}\int_{\R^3}\Re\Big\langle\frac{1}{i}\nabla v_h \mid
\frac{\partial A}{\partial x_k}(z_0+\eps_h x)\bar v_h\Big\rangle dx=
\int_{\R^3}\Re\Big\langle{\textstyle\frac{1}{i}}\nabla v_0 \mid
\frac{\partial A}{\partial x_k}(z_0)\bar v_0\Big\rangle dx.
$$
Since similar considerations apply to the other terms
that appear in~\eqref{ps-prima}, we can therefore pass to the limit
as $h\to\infty$, to find, for each $k$,
\begin{align}
\label{firstcons}
\int_{\R^3}
&\Big\langle\frac{\partial A}{\partial x_k}(z_0)\mid
A(z_0)\Big\rangle|v_0|^2-\Re\Big\langle{\textstyle\frac{1}{i}}\nabla v_0 \mid
\frac{\partial A}{\partial x_k}(z_0)\bar v_0\Big\rangle dx \notag \\
&+\frac{\partial V}{\partial
x_k}(z_0)\int_{\R^3}\frac{|v_0|^2}{2}dx-\frac{\partial K}{\partial
x_k}(z_0)\int_{\R^3}F(|v_0|^2)dx=0.
\end{align}
Now, as proved in Step I, $v_0$ can be represented as
$v_0(x)=e^{i\Upsilon_{z_0}(x)}U_{z_0}(x)$ where
$U_{z_0}:\R^3\to\C$ solves~\eqref{limU}.
Taking into account that
$$
{\textstyle\frac{1}{i}} \nabla {v_0}(x)=e^{i\Upsilon_{z_0}(x)}
A(z_0)U_{z_0}(x)-i e^{i\Upsilon_{z_0}(x)}\nabla U_{z_0}(x),
$$
for every $x\in\R^3$ we obtain
\begin{align*}
\Big\langle\frac{\partial A}{\partial x_k}(z_0) &\mid
A(z_0)\Big\rangle|v_0(x)|^2
-\Re\Big\langle{\textstyle\frac{1}{i}}\nabla v_0(x)\mid
\frac{\partial A}{\partial x_k}(z_0)\bar v_0(x)\Big\rangle \\
&=\Big\langle\frac{\partial A}{\partial x_k}(z_0)\mid
A(z_0)\Big\rangle |U_{z_0}(x)|^2 \\
&\quad-\Re\Big\langle e^{i \Upsilon_{z_0}(x)}A(z_0)U_{z_0}(x)
-i e^{i\Upsilon_{z_0}(x)}\nabla U_{z_0}(x)\mid
\frac{\partial A}{\partial x_k}(z_0) e^{-i\Upsilon_{z_0}(x)}{\bar U}_{z_0}(x)\Big\rangle \\
&=\Big\langle\frac{\partial A}{\partial x_k}(z_0)\mid
A(z_0)\Big\rangle |U_{z_0}(x)|^2 \\
&\quad-\Re\left(\Big\langle\frac{\partial A}{\partial x_k}(z_0)\mid A(z_0)
\Big\rangle |U_{z_0}(x)|^2-i  \Big\langle\frac{\partial A}{\partial x_k}(z_0)
\mid\nabla U_{z_0}(x)\Big\rangle{\bar U}_{z_0}(x)\right) \\
&=\Big\langle\frac{\partial A}{\partial x_k}(z_0)\mid \Re(i{\bar U}_{z_0}(x)
\nabla U_{z_0}(x))\Big\rangle.
\end{align*}
Hence, equation~\eqref{firstcons} can be rephrased as
\begin{align*}
\Big\langle\frac{\partial A}{\partial x_k}(z_0) &\mid\int_{\R^3} \Re(i\bar U_{z_0}
\nabla U_{z_0})dx\Big\rangle \\
&+\frac{\partial V}{\partial
x_k}(z_0)\int_{\R^3}\frac{|U_{z_0}|^2}{2}dx-\frac{\partial K}{\partial
x_k}(z_0)\int_{\R^3}F(|U_{z_0}|^2)dx=0,
\end{align*}
for every $k=1,2,3$, namely,
\begin{align}
\label{finalId}
\Big\langle\frac{\partial A}{\partial w}(z_0) &\mid\int_{\R^3} \Re(i\bar U_{z_0}
\nabla U_{z_0})dx\Big\rangle \notag \\
&+\frac{\partial V}{\partial
w}(z_0)\int_{\R^3}\frac{|U_{z_0}|^2}{2}dx
-\frac{\partial K}{\partial
w}(z_0)\int_{\R^3}F(|U_{z_0}|^2)dx=0,
\end{align}
for every $w\in\R^3$.
\vskip2pt
\noindent
\textbf{Step V.} In this final step, we prove the desired inclusions
stated by the theorem. As a consequence of identity~\eqref{finalId},
in light of the definition of $\Gamma^\pm(z_0;w)$,
we immediately deduce that $z_0\in\mathfrak{S}^*$,
thus proving that $\mathscr{C}\subset\mathfrak{S}^*$.
Let us now assume that $z_0\in{\mathscr E}$. Then
$J_{z_0}(v_0)=\Sigma_c(z_0)=\Sigma_r(z_0)$,
and by virtue of (iii) of Lemma~\ref{pro-sigma1},
we have $U_{z_0}(x)=e^{i\omega}u_{z_0}(x)$
for some $\omega\in\R$, where $u_{z_0}$ is a
real least energy solution to~\eqref{pb:frozen}. Moreover,
by (ii) of Lemma~\ref{pro-sigma1}, we have
$$
\Re\left(i{\bar U}_{z_0}(x)\nabla U_{z_0}(x)\right)=0,
\qquad\text{for a.e.\ $x\in\R^3$.}
$$
Then, in light of Lemma~\ref{pro-sigma0} and~\eqref{finalId}, we obtain
\begin{align*}
\left(\frac{\partial \Sigma_r}{\partial w}\right)^{-}(z_0) &=
\sup_{u\in S_r(z_0)}\Big[\frac{\partial V}{\partial w}(z_0)
\int_{\R^3} \frac{|u|^{2}}{2}dx+
\frac{\partial K}{\partial w}(z_0)
\int_{\R^3} F(|u|^2)dx\Big] \notag \\
&=\sup_{\substack{U=e^{i\omega}u\\
u\in S_r(z_0)}}\Big[\frac{\partial V}{\partial w}(z_0)
\int_{\R^3} \frac{|U|^{2}}{2}dx+
\frac{\partial K}{\partial w}(z_0)
\int_{\R^3} F(|U|^2)dx\Big] \notag \\
&\geq\frac{\partial V}{\partial w}(z_0)
\int_{\R^3} \frac{|U_{z_0}|^{2}}{2}dx+
\frac{\partial K}{\partial w}(z_0)
\int_{\R^3} F(|U_{z_0}|^2)dx=0,
\end{align*}
for every $w\in\R^3$. In a similar fashion, there holds
$$
\left(\frac{\partial \Sigma_r}{\partial w}\right)^{+}(z_0) \leq 0,
$$
for every $w\in\R^3$. In particular, by the definition of $(-\Sigma_r)^0(z_0;w)$, we get
$$
(-\Sigma_r)^0(z_0;w)\geq\left(\frac{\partial (-\Sigma_r)}
{\partial w}\right)^{+}\!\!(z_0)\geq 0,
$$
for $w\in\R^3$. Hence $0\in \partial_C(-\Sigma_r)(z_0)$,
which, in light of Proposition~\ref{oppositesub},
yields $z_0\in\mathfrak{S}$. Finally, if $f(u)$ satisfies~\eqref{fpowernon},
problem~\eqref{pb:frozen} admits a unique real solution $\phi_0$
up to translations (see \cite{chen-lin}).
Taking into account Lemma~\ref{pro-sigma1}, there exists $\omega\in\R$ such that
$v_0=e^{i\Upsilon_{z_0}(x)+i\omega}\phi_0(x)$.
Then, if $z_0\in\mathscr{C}=\mathscr{E}$ (see Proposition~\ref{comparison}), we have
$S_r(z_0)=\{\phi_0\}$, $\Sigma_r$ admits all the directional derivatives and, by the
above inequalities,
$$
\left(\frac{\partial\Sigma_r}{\partial w}\right)^{\pm}\!\!(z_0)=
\frac{\partial\Sigma_r}{\partial w}(z_0)=0,
$$
for $w\in\R^3$. Since up to a multiplicative constant $\Sigma_r$ writes down
explicitly as~\eqref{explicit}, the last assertion readily follows by
a direct computation.
\end{proof}

In light of identity~\eqref{finalId}, we also have the following

\begin{corollary}
Under the assumptions of Theorem~\ref{main}, for every $z_0\in\mathscr{C}$, there exist
constants $\lambda_1,\lambda_2,\lambda_3\in\R$ (possibly zero) and $\gamma_1,\gamma_2
\in\R\setminus\{0\}$ such that
\begin{equation}
\label{alter-cond}
\sum_{j=1}^3\lambda_j\nabla A_j(z_0)
+\gamma_1\nabla V(z_0)+\gamma_2\nabla K(z_0)=0.
\end{equation}
Hence, in general, the location of concentration points
might depend also on the (fixed) external electromagnetic potential $A$.
If $\mathfrak{S}^*=\emptyset$, then~\eqref{problem}
does not admit any sequence of bound-state solutions
concentrating somewhere pointwise.
\end{corollary}

\begin{corollary}
The location of energy-concentration points of a sequence of
bound-state solutions to problem~\eqref{problem} is
independent of the external electromagnetic field $B$
(and there holds $\lambda_j=0$ for all $j=1,2,3$ in~\eqref{alter-cond}).
If $\mathfrak{S}=\emptyset$, then~\eqref{problem}
does not admit any sequence of bound-state solutions
concentrating somewhere energetically.
\end{corollary}

\begin{remark}
\label{rem-gauge}
Despite the fact that both pointwise and energy concentration are
gauge invariant, the necessary condition~\eqref{alter-cond} is not, in general, unless
$\lambda_j=0$ for all $j=1,2,3$. Hence, it seems natural
to conjecture that the answer to Question~\ref{op-quest} is
always affirmative.
\end{remark}

\begin{corollary}
Assume that $f(u)$ is such that, for every $z\in\R^3$, problem~\eqref{pb:frozen} admits a unique
positive radial solution, up to translations. Then, if $z$ is an
energy-concentration point it is a classical critical point of $\Sigma_r$.
\end{corollary}
\noindent
We refer the reader to~\cite[Theorem 2.5 and Theorem 4.2]{chen-lin} for some results
ensuring uniqueness for~\eqref{pb:frozen} under some additional hypothesis on $f(u)$.

\vskip4pt
\noindent
We finish the paper with a simple but interesting
property of the family $\{\mathfrak{S}_p\}_{p\in(1,5)}$.

\begin{observation}
Assume that $f(u)$ satisfies~\eqref{fpowernon} and that
$$
\limsup_{|x|\to\infty}\frac{|\nabla V(x)|}{V(x)}<\infty
\qquad
\text{and}
\qquad
\liminf_{|x|\to\infty}|\nabla K(x)|>0.
$$
We denote by $\operatorname{Crit}(K)$ the set of critical points of $K$, which is
a compact set in light of the above assumption.
Then, it is a simple task to check that
$$
\lim_{p\to 5^-}\operatorname{dist}_{\R^3}(\mathfrak{S}_p,\operatorname{Crit}(K))=0,
$$
that is, if $p$ is close to the critical exponent $5$, the spikes
locate close to $\operatorname{Crit}(K)$.
\end{observation}

\section*{Acknowledgments}
The authors wish to thank Professors Silvia Cingolani and Kazuhiro Kurata for a few comments
about their papers~\cite{cingo} and \cite{ku} respectively.

\end{document}